\documentclass[a4paper,10pt]{amsart}

\usepackage{amsmath,amsfonts,amssymb,amsthm,epsfig}

\def\C{\mathbb{C}}
\def\Z{\mathbb{Z}}
\def\N{\mathbb{N}}

\def\id{\mathbf{1}}

\DeclareMathOperator{\Hom}{Hom}

\DeclareMathOperator{\End}{End}
\DeclareMathOperator{\Aut}{Aut}
\DeclareMathOperator{\inc}{in}
\DeclareMathOperator{\out}{out}
\DeclareMathOperator{\tr}{tr}
\DeclareMathOperator{\Irr}{Irr}
\DeclareMathOperator{\ad}{ad}

\newtheorem{theo}{Theorem}[subsection]
\newtheorem{prop}[theo]{Proposition}
\newtheorem{lem}[theo]{Lemma}
\newtheorem{cor}[theo]{Corollary}

\newtheorem{defin}[theo]{Definition}
\newtheorem{rem}[theo]{Remark}
\numberwithin{equation}{subsection}

\begin{document}
\title[Affine Lie algebras, quiver varieties and statistical mechanics]
{Bases of representations of type $A$ affine Lie algebras via
quiver varieties and statistical mechanics}
\author{Igor B. Frenkel and Alistair Savage}
\date{\today}

\begin{abstract}

We relate two apparently
different bases in the representations of affine Lie algebras of
type $A$: one arising from statistical mechanics, the other from gauge
theory.  We show that the two are governed by the same combinatorics
that also respects the weight space decomposition of the
representations.  In
particular, we are able to give an alternative and much simpler
geometric proof of the main result of \cite{D892} on the construction
of bases of affine Lie algebra representations.  At the same time,
we give a simple parametrization of the irreducible components of
Nakajima quiver varieties associated to infinite and cyclic quivers.
We also define new varieties whose irreducible components are in
one-to-one correspondence with bases of the highest weight representations of
$\widehat{\mathfrak{gl}}_{n+1}$.

\end{abstract}

\maketitle
\noindent

\section*{Introduction}
A remarkable relation between representation theory of affine Lie
algebras and models of statistical mechanics based on the Yang-Baxter
equation has been discovered and intensively studied by E. Date,
M. Jimbo, A. Kuniba, T. Miwa and M. Okado (see \cite{D89,D892} and
references therein).  One of the important findings of the above
authors is that the one-dimensional configuration sums for these
models give rise to characters of integrable highest weight
representations of affine Lie algebras.  This relation yields certain
explicit bases in the representations that admit pure combinatorial
descriptions and imply various identities for the characters.

Another astonishing relation between representation theory of affine
Lie algebras and moduli spaces of solutions of self-dual Yang-Mills
equations has been accomplished by H. Nakajima \cite{N94,N98}, who
observed a profound link between his earlier work with P. Kronheimer
and the results of G. Lusztig \cite{L91,L92}.  At the heart of both
works that preceded the Nakajima discovery are quiver varieties
associated with extended Dynkin diagrams.  Nakajima introduced a
special class of quiver varieties associated with integrable highest
weight representations of affine Lie algebras and obtained a geometric
description of the action.  He also defined certain Lagrangian
subvarieties whose irreducible components yield a geometric basis of
the affine Lie algebra representations.

The central goal of the present paper is to relate the two apparently
different bases in the representations of affine Lie algebras of
type $A$: one arising from statistical mechanics, the other from gauge
theory.  We show that the two are governed by the same combinatorics
that also respects the weight space decomposition of the
representations.  This identification allows
one to give a natural conceptual framework to the intricate
structure of statistical mechanical models and also to make explicit
calculations in a seemingly intractable geometric setting.  In
particular, we are able to give an alternative and much simpler
geometric proof of the main result of \cite{D892} on the construction
of a basis of affine Lie algebra representations.  At the same time,
we give a simple parametrization of the irreducible components of
Nakajima quiver varieties associated to infinite and cyclic quivers.

The comparison of the two very different theories brings some
surprises and suggests interesting new directions.  In particular, the
Young diagrams that are routinely used in representation theory of
type $A$ Lie algebras acquire an explicit geometric meaning:  They
picture precisely representations of the corresponding quivers
satisfying a stability condition for level 1 (see
Figure~\ref{fig:youngdiag} in the text).  On the other hand, the
algebraic constructions of \cite{D892} involve substantially the
highest weight representations of $\widehat{\mathfrak{gl}}_{n+1}$,
which are not directly covered by Nakajima's theory.  We define new
varieties by relaxing the nilpotency condition in the definition of
Nakajima's quiver varieties and show that the irreducible components
of these new varieties are in one-to-one correspondence with bases of the
highest weight representations of $\widehat{\mathfrak{gl}}_{n+1}$.  We
also mention some interesting problems that arise as a result of the
comparision of geometric and algebraic constructions.

We strongly believe that the main results of the current paper reflect a very
general principle that asserts the profound geometric or gauge
theoretic origin of various algebraic and combinatorial structures of
integrable models in statistical mechanics.  The relation of both
subjects to the representation theory of affine Lie algebras is a
necessary prerequisite of this principle.  However we expect much
more; namely that various specific constructions appearing in
integrable models of statistical mechanics that include tensor
products, fusion products, branching rules, Bethe's ansatz and the
Yang-Baxter equation itself reflect certain geometric facts about
Nakajima varieties, Malkin-Nakajima tensor product varieties, various
Lagrangian subvarieties and corresponding gauge theories on
commutative and, possibly, noncommutative spaces.  The present paper
is a small but indicative step toward this vast program.

The paper is organized as follows.  In Section~\ref{sec:lus_def} we recall
the definition of Lusztig's quiver varieties and characterizations of the
irreducible components in types $A_\infty$ and
$A_n^{(1)}$.  We also introduce a version of Lusztig's quiver
varieties for the Lie algebra $\widehat{\mathfrak{gl}}_{n+1}$.
Section~\ref{sec:def_nak} contains the definition of Nakajima's quiver
varieties and the Lie algebra action on a suitable space of
constructible functions on these varieties is given in
Section~\ref{sec:structure}.  In Section~\ref{sec:level1} we give an
enumeration of the irreducible components of the quiver varieties for
level 1 in terms of Young diagrams.  We also identify the geometric
action of the type $A_\infty$ Lie algebra in the basis enumerated by
Young diagrams.  In Section~\ref{sec:arblevel} we extend the
enumeration of the irreducible components of the quiver varieties to
arbitrary level and we establish a match with the indexing of bases of
the corresponding representations coming from statistical mechanics.
Finally, in Section~\ref{sec:comparison}, we compare the weight
structure of the bases resulting from quiver varieties and the path
realizations of statistical mechanics and make certain of their complete
coincidence.

The research of the first author was supported in part by the National
Science Foundation (NSF).  The research of the second author was supported
in part by the Natural Sciences and Engineering Research Council of
Canada (NSERC).


\section{Lusztig's Quiver Varieties}
\label{sec:lus_def}
\setcounter{theo}{0}

In this section, we will recount the explicit description given in
\cite{L91} of the
irreducible components of Lusztig's quiver variety in the case of
types $A_{\infty}$ and
$A_n^{(1)}$.  See this reference for details, including proofs.

Let $I$ be a set of vertices of the Dynkin graph of a Kac-Moody Lie
algebra $\mathfrak{g}$ and let $H$ be the set of pairs consisting of
an edge together with an orientation of it.  For $h \in H$, let
$\inc(h)$ (resp. $\out(h)$) be the incoming (resp. outgoing) vertex of
$h$.  We define the involution $\bar{\ }: H \to H$ to be the function
which takes $h \in H$ to the element of $H$ consisting of the the same
edge with opposite orientation.  An \emph{orientation} of
our graph is a choice of a subset $\Omega \subset H$ such that $\Omega
\cup \bar{\Omega} = H$ and $\Omega \cap \bar{\Omega} = \emptyset$.

Let $\mathcal{V}$ be the category of
finite-dimensional $I$-graded vector spaces $\mathbf{V} = \oplus_{i
  \in I} \mathbf{V}_i$ over $\C$ with morphisms being linear maps
respecting the grading.  Then $\mathbf{V} \in \mathcal{V}$ shall denote
that $\mathbf{V}$ is an object of $\mathcal{V}$.  The dimension of
$\mathbf{V} \in \mathcal{V}$ is given by $\mathbf{v} = \dim \mathbf{V}
= (\dim \mathbf{V}_0, \dots,
\dim \mathbf{V}_n)$.  We 
identify this dimension with the element $(\dim \mathbf{V}_0)\alpha_0
+ \dots + (\dim
\mathbf{V}_n) \alpha_n$ of the root lattice of $\mathfrak{g}$.
Here the $\alpha_i$ are the simple roots corresponding to
the vertices of our quiver (graph with orientation), whose underlying graph is
the Dynkin graph of $\mathfrak{g}$.

Given $\mathbf{V} \in \mathcal{V}$, let
\[
\mathbf{E_V} = \bigoplus_{h \in H} \Hom (\mathbf{V}_{\out(h)},
\mathbf{V}_{\inc(h)}).
\]
For any subset $H'$ of $H$, let $\mathbf{E}_{\mathbf{V}, H'}$ be the
subspace of $\mathbf{E_V}$ consisting of all vectors $x = (x_h)$ such
that $x_h=0$ whenever $h \not\in H'$.  The algebraic group $G_\mathbf{V}
= \prod_i \Aut(\mathbf{V}_i)$ acts on $\mathbf{E_V}$ and
$\mathbf{E}_{\mathbf{V}, H'}$ by
\[
(g,x) = ((g_i), (x_h)) \mapsto gx = (x'_h) = (g_{\inc(h)} x_h g_{\out(h)}^{-1}).
\]

Define the function $\varepsilon : H \to \{-1,1\}$ by $\varepsilon (h) =
1$ for all $h
\in \Omega$ and $\varepsilon(h) = -1$ for all $h \in {\bar{\Omega}}$.  Let
$\mathbf{V} \in \mathcal{V}$.  The Lie algebra of $G_\mathbf{V}$ is
$\mathbf{gl_V} = \prod_i \End(\mathbf{V}_i)$ and it acts on $\mathbf{E_V}$ by
\[
(a,x) = ((a_i), (x_h)) \mapsto [a,x] = (x'_h) = (a_{\inc(h)}x_h - x_h
a_{\out(h)}).
\]
Let $\left<\cdot,\cdot\right>$ be the nondegenerate,
$G_\mathbf{V}$-invariant, symplectic form on
$\mathbf{E_V}$ with values in $\C$ defined by
\[
\left<x,y\right> = \sum_{h \in H} \varepsilon(h) \tr (x_h y_{\bar{h}}).
\]
Note that $\mathbf{E_V}$ can be considered as the cotangent space of
$\mathbf{E}_{\mathbf{V}, \Omega}$ under this form.

The moment map associated to the $G_{\mathbf{v}}$-action on the
symplectic vector space $\mathbf{E_V}$ is the map $\psi : \mathbf{E_V}
\to \mathbf{gl_V}$ with $i$-component $\psi_i : \mathbf{E_V} \to \End
\mathbf{V}_i$ given by
\[
\psi_i(x) = \sum_{h \in H,\, \inc(h)=i} \varepsilon(h) x_h x_{\bar{h}} .
\]

\begin{defin}[\cite{L91}]
\label{def:nilpotent}
An element $x \in \mathbf{E_V}$ is said to be \emph{nilpotent} if
there exists an $N \ge 1$ such that for any sequence $h'_1, h'_2,
\dots, h'_N$ in $H$ satisfying $\out (h'_1) = \inc (h'_2)$, $\out (h'_2) =
\inc (h'_3)$, \dots, $\out (h'_{N-1}) = \inc (h'_N)$, the composition
$x_{h'_1} x_{h'_2} \dots x_{h'_N} : \mathbf{V}_{\out (h'_N)} \to
  \mathbf{V}_{\inc (h'_1)}$ is zero.
\end{defin}

\begin{defin}[\cite{L91}] Let $\Lambda_\mathbf{V}$ be the set of all
  nilpotent elements $x \in \mathbf{E_V}$ such that $\psi_i(x) = 0$
  for all $i \in I$.
\end{defin}


\subsection{Type $A_\infty$}
Let $\mathfrak{g}$ be the simple Lie algebra of type $A_\infty$. Let
$I=\Z$ be the set of vertices of a graph with the set of oriented
edges given by
\begin{gather*}
H=\{
i \to j \ |\ i,j \in I,\ i-j=1\} \cup
\{i \leftarrow j \ |\ i,j \in I, i-j=1\}.
\end{gather*}

We define the involution $\bar{\ } : H \to H$ as the function that
interchanges $i \to j$ and $i \leftarrow j$.  For $h=(i \to j)$, we
set $\out(h) = i$ and $\inc(h) = j$ and for $h=(i \leftarrow j)$, we
set $\out(h) = j$ and $\inc(h) = i$.  Let $\Omega$ be the subset of
$H$ consisting of the arrows $i \to j$.

\begin{prop}[\cite{L91}]
\label{prop:irrcomp:Ainfty}
The irreducible components of $\Lambda_\mathbf{V}$ are the closures of
the conormal bundles of the various $G_\mathbf{V}$-orbits in
$\mathbf{E}_{\mathbf{V}, \Omega}$.
\end{prop}
\begin{proof}
The case where $\mathfrak{g}$ is of type $A_n$ is proven in
\cite{L91}.  The $A_\infty$ case follows by passing to the direct
limit.
\end{proof}

For two integers $k' \le k$, define $\mathbf{V}_\infty(k',k) \in \mathcal{V}$
to be the vector space with basis $\{ e_r\ |\ k' \le r \le k\}$.  We
require that $e_r$ has degree $r \in I$.
Let $x_\infty(k', k) \in \mathbf{E}_{\mathbf{V}_\infty(k',k), \Omega}$
be defined by
$x_\infty(k',k) : e_r \mapsto e_{r-1}$ for $k' \le r \le k$, where $e_{k'-1}
= 0$.  It is clear that $(\mathbf{V}_\infty(k',k), x_\infty(k',k))$ is an
indecomposable representation of our quiver.
Conversely, any indecomposable finite-dimensional representation
$(\mathbf{V},x)$ of our quiver is isomorphic to
some $(\mathbf{V}_\infty(k',k),x_\infty(k',k))$.

Let $Z^\infty$ be the set of all pairs $(k' \le k)$ of integers and
let $\tilde Z^\infty$ be
the set of all functions $Z^\infty \to \N$ with finite support.

It is easy to see that for $\mathbf{V} \in \mathcal{V}$, the set of
$G_\mathbf{V}$-orbits in
$\mathbf{E}_{\mathbf{V}, \Omega}$ is naturally indexed by the subset
$\tilde Z^\infty_\mathbf{V}$ of $\tilde Z^\infty$ consisting of those
$f \in \tilde Z^\infty$ such that
\[
\sum_{k' \le i \le k} f(k',k) = \dim \mathbf{V}_i
\]
for all $i \in I$.  Here the sum is over all $k' \le k$ such that $k'
\le i \le k$.  Corresponding to a given $f$ is the orbit consisting of
all representations isomorphic to a sum of the indecomposable
representations $x_\infty(k',k)$, each occuring with multiplicity $f(k',k)$.
Denote by
$\mathcal{O}_f$ the $G_\mathbf{V}$-orbit corresponding to $f \in
\tilde Z^\infty_\mathbf{V}$.  Let $\mathcal{C}_f$ be the
conormal bundle to $\mathcal{O}_f$ and let
$\bar{\mathcal{C}}_f$ be its closure.  We then have the following
proposition.

\begin{prop}
The map $f \to \bar{\mathcal{C}}_f$ is a one-to-one correspondence
between the set
$\tilde Z^\infty_\mathbf{V}$ and the set of irreducible
components of $\Lambda_\mathbf{V}$.
\end{prop}
\begin{proof}
This follows immediately from Proposition~\ref{prop:irrcomp:Ainfty}.
\end{proof}


\subsection{Type $A_n^{(1)}$}
Let $\mathfrak{g}$ be the affine Lie algebra of type $A_n^{(1)}$,
that is, the Lie algebra generated by the set of elements $E_k$,
$F_k$, $H_k$ ($k=0,1,\dots,n$) and $d$ satisfying the following
relations:
\begin{gather*}
[E_k,F_l] = \delta_{kl} H_k,\qquad [H_k,E_l]=a_{kl}E_l,\qquad [H_k,F_l]=-a_{kl}F_l,
\\
[d,E_k] = \delta_{k0}E_k,\qquad [d,F_k]=-\delta_{k0}F_k, \\
(\ad E_k)^{1-a_{kl}}E_l = 0,\qquad (\ad F_k)^{1-a_{kl}}F_l = 0 \quad
\text{ for }
k \ne l.
\end{gather*}
Here
\[
a_{kl} = 2\delta(k,l) - \delta(k,l+1) - \delta(k,l-1),
\]
where $\delta(k,l) = 1$ if $k \equiv l$ mod $(n+1)$ and is equal to zero
otherwise.

Let $I=\Z/(n+1)\Z$ be the set of vertices of a graph with the set of
oriented edges given by
\begin{gather*}
H=\{
i \to j \ |\ i,j \in I,\ i-j=1\} \cup
\{i \leftarrow j \ |\ i,j \in I, i-j=1\}.
\end{gather*}

For two integers $k' \le k$, define $\mathbf{V}(k',k) \in \mathcal{V}$
to be the vector space with basis $\{ e_r\ |\ k' \le r \le k\}$.  We
require that $e_r$ has degree $i \in I$ where $r \equiv i$ (mod $n+1$).
Let $x(k', k) \in \mathbf{E}_{\mathbf{V}(k',k), \Omega}$ be defined by
$x(k',k) : e_r \mapsto e_{r-1}$ for $k' \le r \le k$, where $e_{k'-1}
= 0$.  It is clear that $(\mathbf{V}(k',k), x(k',k))$ is an
indecomposable representation of our quiver and that $x(k',k)$ is nilpotent.
Also, the isomorphism class of this representation does not change when
$k'$ and $k$ are simultaneously translated by a multiple of $n+1$.
Conversely, any indecomposable finite-dimensional representation
$(\mathbf{V},x)$ of our quiver, with $x$ nilpotent, is isomorphic to
some $(\mathbf{V}(k',k),x(k',k))$ where $k'$ and $k$ are uniquely
defined up to a simultaneous translation by a multiple of $n+1$.

Let $Z$ be the set of all pairs $(k' \le k)$ of integers defined up to
simultaneous translation by a multiple of $n+1$ and let $\tilde Z$ be
the set of all functions $Z \to \N$ with finite support.

It is easy to see that for $\mathbf{V} \in \mathcal{V}$, the set of
$G_\mathbf{V}$-orbits on the set of nilpotent elements in
$\mathbf{E}_{\mathbf{V}, \Omega}$ is naturally indexed by the subset
$\tilde Z_\mathbf{V}$ of $\tilde Z$ consisting of those $f \in \tilde
Z$ such that
\[
\sum_{k' \le k} f(k',k) \#\{r\ |\ k' \le r \le k,\ r \equiv i \text{
  (mod $n+1$)}\} = \dim \mathbf{V}_i
\]
for all $i \in I$.  Here the sum is taken over all $k' \le k$ up to
simultaneous translation by a multiple of $n+1$.  Corresponding to a
given $f$ is the orbit consisting of
all representations isomorphic to a sum of the indecomposable
representations $x(k',k)$, each occuring with multiplicity
$f(k',k)$. Denote by
$\mathcal{O}_f$ the $G_\mathbf{V}$-orbit corresponding to $f \in
\tilde Z_\mathbf{V}$.

We say that $f \in \tilde Z_\mathbf{V}$ is \emph{aperiodic} if for any
$k' \le k$, not all $f(k',k)$, $f(k'+1, k+1)$, \dots, $f(k'+n, k+n)$
are greater than zero.  For any $f \in \tilde Z_\mathbf{V}$, let
$\mathcal{C}_f$ be
the conormal bundle of $\mathcal{O}_f$ and let $\bar{\mathcal{C}}_f$ be
its closure.

\begin{prop}[{\cite[15.5]{L91}}]
Let $f \in \tilde Z_{\mathbf{V}}$.  The following two conditions are
equivalent.
\begin{enumerate}
\item $\mathcal{C}_f$ consists entirely of
  nilpotent elements.
\item $f$ is aperiodic.
\end{enumerate}
\end{prop}

\begin{prop}[{\cite[15.6]{L91}}]
The map $f \to \bar{\mathcal{C}}_f$ is a 1-1 correspondence between the set of
aperiodic elements in $\tilde Z_\mathbf{V}$ and the set of irreducible
components of $\Lambda_\mathbf{V}$.
\end{prop}

\begin{prop}[{\cite[12.8]{L91}}]
\label{prop:orbit_char}
Let $x' \in \mathbf{E}_{\mathbf{V},\Omega}$ and $x'' \in
\mathbf{E}_{\mathbf{V},\bar{\Omega}}$.  Then $\psi_i(x' + x'')=0$ for
all $i \in I$ if and only if $x''$ is orthogonal with respect to
$\left<\ ,\ \right>$ to the tangent space to the $G_\mathbf{V}$-orbit
  of $x'$, regarded as a vector subspace of $\mathbf{E}_{\mathbf{V},
    \Omega}$.
\end{prop}


\subsection{$\widehat{\mathfrak{gl}}_{n+1}$ Case}
\label{sec:rem_gl_lus}

Since $\widehat{\mathfrak{gl}}_{n+1}$ is not a Kac-Moody algebra in a
strict sense, this case is not covered by Lusztig's theory and
requires certain modifications.  We preserve the notation of the
previous subsection.

\begin{defin}
Let $\tilde \Lambda_\mathbf{V}$ be the set of all elements $x=x'+x''$,
where $x'
\in E_{\mathbf{V},\Omega}$ and $x'' \in E_{\mathbf{V},\bar{\Omega}}$,
such that $x'$ is nilpotent and $\psi_i(x)=0$ for all $i
\in I$.
\end{defin}

For any $f \in \tilde Z_\mathbf{V}$, we denote by $\mathcal{O}_f$ the
corresponding $G_\mathbf{V}$-orbit and by $\mathcal{C}_f$ its conormal
bundle.

\begin{prop}
Let $f \in \tilde Z_\mathbf{V}$.  Then
\begin{enumerate}
\item $\mathcal{C}_f$ consists entirely of elements of $\tilde
  \Lambda_\mathbf{V}$, and
\item $\tilde \Lambda_\mathbf{V}$ is the union of ${\bar{\mathcal{C}}}_f$ for
  all $f \in \tilde Z_\mathbf{V}$.
\end{enumerate}
\end{prop}

\begin{proof}
This follows from Proposition~\ref{prop:orbit_char}.
\end{proof}

\begin{prop}
The map $f \to \bar{\mathcal{C}}_f$ is a 1-1 correspondence between the set
$\tilde Z_\mathbf{V}$ and the set of irreducible
components of $\Lambda_\mathbf{V}$.
\end{prop}

\begin{proof}
This follows easily since the conormal bundles $\mathcal{C}_f$ are
irreducible of the same dimension. 
\end{proof}


\section{Nakajima's Quiver Varieties}
\label{sec:def_nak}
\setcounter{theo}{0}

We introduce here a description of the quiver varieties first
presented in \cite{N94} in the case of types $A_\infty$ and $A_n^{(1)}$.

\begin{defin}[\cite{N94}]
\label{def:lambda}
For $\mathbf{v}, \mathbf{w} \in \Z_{\ge 0}^I$, choose $I$-graded
vector spaces $\mathbf{V}$ and $\mathbf{W}$ of graded dimensions
$\mathbf{v}$ and
$\mathbf{w}$ respectively.  Then define
\[
\Lambda \equiv \Lambda(\mathbf{v},\mathbf{w}) =
\Lambda_\mathbf{V} \times \sum_{i \in I} \Hom (\mathbf{V}_i, \mathbf{W}_i).
\]
\end{defin}

Now, suppose that $\mathbf{S}$ is an $I$-graded subspace of $\mathbf{V}$.
For $x \in
\Lambda_\mathbf{V}$ we say that $\mathbf{S}$ is
\emph{$x$-stable} if $x(\mathbf{S}) \subset \mathbf{S}$.

\begin{defin}[\cite{N94}]
\label{def:lambda-stable}
Let $\Lambda^{\text{st}} = \Lambda(\mathbf{v},\mathbf{w})^{\text{st}}$ be
the set of all $(x, j) \in
\Lambda(\mathbf{v},\mathbf{w})$ satisfying the following condition:  If
$\mathbf{S}=(\mathbf{S}_i)$ with $\mathbf{S}_i \subset \mathbf{V}_i$ is
$x$-stable and $j_i(\mathbf{S}_i) = 0$ for
$i \in I$, then $\mathbf{S}_i = 0$ for $i \in I$.
\end{defin}

The group $G_\mathbf{V}$ acts on $\Lambda(\mathbf{v},\mathbf{w})$ via
\[
(g,(x,j)) = ((g_i), ((x_h), (j_i))) \mapsto ((g_{\inc (h)} x_h
g_{\out (h)}^{-1}), j_i g_i^{-1}).
\]
and the stabilizer of any point of
$\Lambda(\mathbf{v},\mathbf{w})^{\text{st}}$ in $G_{\mathbf{V}}$ is trivial
  (see \cite[Lemma~3.10]{N98}).  We then make the following definition.
\begin{defin}[\cite{N94}]
\label{def:L}
Let $\mathcal{L} \equiv \mathcal{L}(\mathbf{v},\mathbf{w}) =
\Lambda(\mathbf{v},\mathbf{w})^{\text{st}} / G_{\mathbf{V}}$.
\end{defin}

Let $\Irr \mathcal{L}(\mathbf{v},\mathbf{w})$ (resp. $\Irr
\Lambda(\mathbf{v},\mathbf{w})$) be the set of irreducible components of
$\mathcal{L}(\mathbf{v},\mathbf{w})$ (resp. $\Lambda(\mathbf{v},\mathbf{w})$).
Then $\Irr \mathcal{L}(\mathbf{v},\mathbf{w})$ can be identified with
\[
\{ Y \in \Irr \Lambda(\mathbf{v},\mathbf{w})\, |\, Y \cap
\Lambda(\mathbf{v},\mathbf{w})^{\text{st}} \ne \emptyset \}.
\]
Specifically, the irreducible components of $\Irr
\mathcal{L}(\mathbf{v},\mathbf{w})$ are precisely those
\[
X_f \stackrel{\text{def}}{=} \left( \left( \bar{\mathcal{C}}_f \times
    \sum_{i \in I} \Hom (\mathbf{V}_i, \mathbf{W}_i) \right) \cap
\Lambda(\mathbf{v},\mathbf{w})^{\text{st}} \right) / G_\mathbf{V}
\]
which are nonempty.

The following will be used in the sequel.
\begin{lem}
\label{lem:irrcomp}
One has
\[
\Lambda^{\text{st}} = \{ x \in \Lambda\, |\, \ker x_{i
\to i-1} \cap \ker x_{i+1 \leftarrow i} \cap \ker j_i = 0 \ \forall i\}
\].
\end{lem}
\begin{proof}
Since each $\ker x_{i \to i-1} \cap \ker x_{i+1 \leftarrow i}$ is
$x$-stable, the left hand side is obviously
contained in the right hand side.  Now suppose $x$ is an element of the
right hand side.  Let $\mathbf{S} = (\mathbf{S}_i)$ with $\mathbf{S}_i
\subset \mathbf{V}_i$ be $x$-stable
and $j_i (\mathbf{S}_i) = 0$ for $i \in I$.  Assume that $\mathbf{S} \ne 0$.
Since all elements of $\Lambda$ are nilpotent, we can
find a minimal value of $N$ such that the condition in
Definition~\ref{def:nilpotent} is satisfied.  Then we can find a $v
\in \mathbf{S}_i$ for some $i$ and a sequence
$h_1', h_2', \dots, h_{N-1}'$ (empty if $N=1$) in $H$ such that $\out
(h_1') = \inc
(h_2')$, $\out (h_2') = \inc (h_3')$, \dots, $\out (h_{N-2}') = \inc
(h_{N-1}')$ and $v' = x_{h_1'} x_{h_2'} \dots x_{h_{N-1}'} (v) \ne 0$.
Now, $v' \in \mathbf{S}_{i'}$ for some $i' \in I$ by the stability of
$\mathbf{S}$ (hence
$j_{i'}(v') = 0$) and $v' \in \ker x_{i' \to i'-1} \cap \ker x_{i'+1
  \to i'}$ by our choice of $N$.  This contradicts the fact
that $x$ is an element of the right hand side.
\end{proof}

In the case of $\widehat{\mathfrak{gl}}_{n+1}$, we define the
varieties $\tilde{\Lambda}(\mathbf{v},\mathbf{w})$,
$\tilde{\Lambda}(\mathbf{v},\mathbf{w})^{st}$ and
$\tilde{\mathcal{L}}(\mathbf{v}, \mathbf{w})$ by
replacing $\Lambda_{\mathbf{V}}$ by $\tilde{\Lambda}_{\mathbf{V}}$ in the
above.


\section{The Lie Algebra Action}
\label{sec:structure}
\setcounter{theo}{0}

We summarize here some results from \cite{N94} that will be needed in
the sequel.  See this reference for more details, including proofs.
We keep the notation of Sections~\ref{sec:lus_def}
and~\ref{sec:def_nak} (with $\mathfrak{g}$ arbitrary).

Let $\mathbf{w, v, v', v''} \in \Z_{\ge 0}^I$ be such that $\mathbf{v}
= \mathbf{v'} + \mathbf{v''}$.  Consider the maps
\begin{equation}
\label{eq:diag_action}
\Lambda(\mathbf{v}'',\mathbf{0}) \times \Lambda(\mathbf{v}',\mathbf{w})
\stackrel{p_1}{\leftarrow} \mathbf{\tilde F (v,w;v'')}
\stackrel{p_2}{\rightarrow} \mathbf{F(v,w;v'')}
\stackrel{p_3}{\rightarrow} \Lambda(\mathbf{v},\mathbf{w}), 
\end{equation}
where the notation is as follows.  A point of $\mathbf{F(v,w;v'')}$ is
a point $(x,j) \in \Lambda(\mathbf{v},\mathbf{w})$ together with an $I$-graded,
$x$-stable
subspace $\mathbf{S}$ of $\mathbf{V}$ such that $\dim \mathbf{S} =
\mathbf{v'} = \mathbf{v} - \mathbf{v''}$.  A point of $\mathbf{\tilde
  F (v,w;v'')}$ is a point $(x,j,\mathbf{S})$ of $\mathbf{F(v,w;v'')}$
together with a collection of isomorphisms $R'_i : \mathbf{V}'_i \cong
\mathbf{S}_i$ and $R''_i : \mathbf{V}''_i \cong \mathbf{V}_i /
\mathbf{S}_i$ for each $i \in I$.  Then we define $p_2(x,j,\mathbf{S},
R',R'') = (x,j,\mathbf{S})$, $p_3(x,j,\mathbf{S}) = (x,j)$ and
$p_1(x,j,\mathbf{S},R',R'') = (x'',x',j')$ where $x'', x', j'$ are
determined by
\begin{align*}
R'_{\inc(h)} x'_h &= x_h R'_{\out(h)} : \mathbf{V}'_{\out(h)} \to
\mathbf{S}_{\inc(h)}, \\
j'_i &= j_i R'_i : \mathbf{V}'_i \to \mathbf{W}_i \\
R''_{\inc(h)} x''_h &= x_h R''_{\out(h)} : \mathbf{V}''_{\out(h)} \to
\mathbf{V}_{\inc(h)} / \mathbf{S}_{\inc(h)}.
\end{align*}
It follows that $x'$ and $x''$ are nilpotent.

\begin{lem}[{\cite[Lemma 10.3]{N94}}]
One has
\[
(p_3 \circ p_2)^{-1} (\Lambda(\mathbf{v},\mathbf{w})^{\text{st}}) \subset
p_1^{-1} (\Lambda(\mathbf{v}'',\mathbf{0}) \times
\Lambda(\mathbf{v}',\mathbf{w})^{\text{st}}).
\]
\end{lem}

Thus, we can restrict \eqref{eq:diag_action} to
$\Lambda^{\text{st}}$, forget the
$\Lambda(\mathbf{v}'',\mathbf{0})$-factor and consider the quotient by
$G_\mathbf{V}$, $G_\mathbf{V'}$.  This yields the diagram
\begin{equation}
\label{eq:diag_action_mod}
\mathcal{L}(\mathbf{v'}, \mathbf{w}) \stackrel{\pi_1}{\leftarrow}
\mathcal{F}(\mathbf{v}, \mathbf{w}; \mathbf{v} - \mathbf{v'})
\stackrel{\pi_2}{\rightarrow} \mathcal{L}(\mathbf{v}, \mathbf{w}),
\end{equation}
where
\[
\mathcal{F}(\mathbf{v}, \mathbf{w}, \mathbf{v} - \mathbf{v'})
\stackrel{\text{def}}{=} \{ (x,j,\mathbf{S}) \in \mathbf{F(v,w;v-v')}\,
  |\, (x,j) \in \Lambda(\mathbf{v},\mathbf{w})^{\text{st}} \} / G_\mathbf{V}.
\]

Let $M(\mathcal{L}(\mathbf{v}, \mathbf{w}))$ be the vector space of
all constructible functions on $\mathcal{L}(\mathbf{v}, \mathbf{w})$.
For a subvariety $Y$ of a variety $A$, let $\mathbf{1}_Y$ denote the
function on $A$ which takes the value 1 on $Y$ and 0 elsewhere.  Let
$\chi (Y)$ denote the Euler characteristic of the algebraic variety
$Y$.  Then for a map $\pi$ between algebraic varieties $A$ and $B$, let
$\pi_!$ denote the map between the abelian groups of constructible
functions on $A$ and $B$ given by
\[
\pi_! (\mathbf{1}_Y)(y) = \chi (\pi^{-1}(y) \cap Y),\ Y \subset A
\]
and let $\pi^*$ be the pullback map from functions on $B$ to functions
on $A$ acting as $\pi^* f(y) = f(\pi(y))$.
Then define
\begin{align*}
&H_i : M(\mathcal{L}(\mathbf{v}, \mathbf{w})) \to
M(\mathcal{L}(\mathbf{v}, \mathbf{w})); \quad H_i f = u_i f, \\
&E_i : M(\mathcal{L}(\mathbf{v}, \mathbf{w})) \to
M(\mathcal{L}(\mathbf{v} - \mathbf{e}^i, \mathbf{w})); \quad E_i f =
(\pi_1)_! (\pi_2^* f), \\
&F_i : M(\mathcal{L}(\mathbf{v} - \mathbf{e}^i, \mathbf{w})) \to
M(\mathcal{L}(\mathbf{v}, \mathbf{w})); \quad F_i g = (\pi_2)_!
(\pi_1^* g).
\end{align*}
Here
\[
\mathbf{u} = {^t(u_0, \dots, u_n)} = \mathbf{w} - C \mathbf{v}
\]
where $C$ is the Cartan matrix of $\mathfrak{g}$ and we are using
diagram~\eqref{eq:diag_action_mod} with $\mathbf{v}' = \mathbf{v} -
\mathbf{e}^i$ where $\mathbf{e}^i$ is the vector whose components are
given by $\mathbf{e}^i_{i'} = \delta_{ii'}$.

Now let $\varphi$ be the constant function on $\mathcal{L}(\mathbf{0},
\mathbf{w})$ with value 1.  Let $L(\mathbf{w})$ be the vector space of
functions generated by acting on $\varphi$ with all possible combinations of
the operators $F_i$.  Then let $L(\mathbf{v},\mathbf{w}) =
M(\mathcal{L}(\mathbf{v}, \mathbf{w})) \cap L(\mathbf{w})$.

\begin{prop}[{\cite[Thm 10.14]{N94}}]
The operators $E_i$, $F_i$, $H_i$ on $L(\mathbf{w})$ provide the
structure of the
irreducible highest weight integrable representation of
$\mathfrak{g}$ with highest weight $\mathbf{w}$.  Each
summand of the decomposition $L(\mathbf{w}) = \bigoplus_\mathbf{v}
L(\mathbf{v}, \mathbf{w})$ is a weight space with weight
$\mathbf{w} - C \mathbf{v}$.
\end{prop}

Let $X \in \Irr \mathcal{L}(\mathbf{v}, \mathbf{w})$ and define a
linear map $T_X : L(\mathbf{v}, \mathbf{w}) \to \C$ as in
\cite[3.8]{L92}.  The map $T_X$ associates to a constructible function $f \in
L(\mathbf{v}, \mathbf{w})$ the (constant) value of $f$ on a suitable
open dense subset of $X$.  The fact that $L(\mathbf{v}, \mathbf{w})$
is finite-dimensional allows us to take such an open set on which
\emph{any} $f \in L(\mathbf{v}, \mathbf{w})$ is constant.  So we have
a linear map
\[
\Phi : L(\mathbf{v}, \mathbf{w}) \to \C^{\Irr \mathcal{L}(\mathbf{v},
  \mathbf{w})}.
\]
The following proposition is proved in \cite[4.16]{L92} (slightly
generalized in
\cite[Proposition 10.15]{N94}).

\begin{prop}
\label{prop:func_irrcomp_isom}
The map $\Phi$ is an isomorphism; for any $X \in \Irr \mathcal{L}(\mathbf{v},
\mathbf{w})$, there is a unique function $g_X \in L(\mathbf{v},
\mathbf{w})$ such that for some open dense subset $O$ of $X$ we have
$g_X|_O = 1$ and for some closed $G_\mathbf{V}$-invariant
subset $K \subset \mathcal{L}(\mathbf{v}, \mathbf{w})$ of dimension $<
\dim \mathcal{L}(\mathbf{v}, \mathbf{w})$ we have $g_X=0$ outside $X
\cup K$.  The functions $g_X$ for $X \in \Irr \Lambda(\mathbf{v},\mathbf{w})$
form a basis of $L(\mathbf{v},\mathbf{w})$.
\end{prop}


\section{Level One Representations}
\label{sec:level1}
\setcounter{theo}{0}

We now seek to describe the irreducible components of Nakajima's
quiver variety.  By the comment made in
Section~\ref{sec:def_nak}, it
suffices to determine which irreducible components of
$\Lambda(\mathbf{v},\mathbf{w})$ are not killed by the stability
condition.  By
Definition~\ref{def:lambda} and Lemma~\ref{lem:irrcomp}, these are
precisely those irreducible components which contain points $x$ such
that
\begin{equation}
\label{eq:killcond}
\dim \left( \ker x_{i \to i-1} \cap \ker x_{i+1 \leftarrow i} \right)
  \le \mathbf{w}_i \ \forall i.
\end{equation}

We first consider the basic representation of
highest weight $\Lambda_0$ where $\Lambda_0(\alpha_i) = \delta_{0i}$.
This corresponds to $\mathbf{w} = \mathbf{w}^0$, the vector with zero
component 1 and all other components equal to zero.


\subsection{Type $A_\infty$}
\label{sec:Ainflevel1}

Consider the case where $\mathfrak{g}$ is of type $A_\infty$.
Let $\mathcal{Y}$ be the set of all Young diagrams, that is,
the set of all weakly decreasing sequences $[l_1, \dots,
l_s]$ of non-negative integers ($l_j=0$ for
$j>s$).  For $Y=[l_1, \dots, l_s] \in \mathcal{Y}$, let $A_Y$ be the
set $\{(1-i,l_i-i)\ |\ 1 \le i \le s\}$.

\begin{theo}
\label{thm:irrcomp_lev1_infty}
The irreducible components of $\mathcal{L}(\mathbf{v},\mathbf{w}^0)$
are precisely those $X_f$ where $f \in \tilde Z^\infty_\mathbf{V}$ such that
\[
\{(k',k)\ |\ f(k',k)=1\} = A_Y
\]
for some $Y \in \mathcal{Y}$ and $f(k',k)=0$ for $(k',k) \not \in
A_Y$.
Denote the component corresponding to such an $f$ by $X_Y$.  Thus,
$Y \leftrightarrow X_Y$ is a natural 1-1 correspondence between the set
$\mathcal{Y}$ and the
irreducible components of $\cup_\mathbf{v} \mathcal{L}(\mathbf{v},
\mathbf{w}^0)$.
\end{theo}

\begin{proof}
Consider the two representations
$(\mathbf{V}_\infty(k_1',k_1),x_\infty(k_1',k_1))$
and $(\mathbf{V}_\infty(k_2',k_2),x_\infty(k_2',k_2))$ of our oriented graph as
described in Section~\ref{sec:lus_def} where the basis of
$\mathbf{V}_\infty(k_i',k_i)$ is $\{e^i_r\ |\ k_i' \le r \le k_i\}$.  Let
$W$ be the conormal bundle to the $G_\mathbf{V}$-orbit through the point
\[
x_\Omega = (x_h)_{h \in \Omega} = x_\infty(k_1',k_1) \oplus
  x_\infty(k_2',k_2) \in
  \mathbf{E}_{\mathbf{V}_\infty(k_1',k_1)
  \oplus \mathbf{V}_\infty(k_2',k_2), \Omega}.
\]
By Proposition~\ref{prop:orbit_char}, $x = (x_\Omega, x_{\bar \Omega})
= (x_h)_{h \in H}$ is in $W$ if
and only if
\[
x_{i+1 \to i} x_{i+1 \leftarrow i} = x_{i \leftarrow i-1} x_{i \to i-1}
\]
for all $i$.

Let $e^i_r=0$ for $r< k_i'$ or $r > k_i$.  Now, $x_{r+1 \leftarrow
  r}(e^2_r) = c_r e^1_{r+1}$ for some $c_r \in \C$ since $x_{r+1
  \leftarrow r}(e^2_r)$ can have no $e^2_{r+1}$-component by
  nilpotency.  Suppose that $k_1' \le r+1 \le k_1$ and $c_r
  \ne 0$ (that is, $x_{r+1 \leftarrow
  r}(e^2_r) \ne 0$) .  Then if $r+1 > k_1'$,
\[
x_{r \leftarrow r-1} (e^2_{r-1}) =
x_{r \leftarrow r-1} x_{r \to r-1} (e^2_r) = x_{r+1 \to r} x_{r+1
  \leftarrow r} (e^2_r) = c_r e^1_r \ne 0.
\]
In particular, $e^2_{r-1} \ne 0$ and so $r-1 \ge k_2'$.  Continuing
in this manner, we see that $x_{k_1' \leftarrow k_1'-1} (e^2_{k_1'-1})
\ne 0$ and thus $k_2' < k_1'$.

Now, if $r+1 \le k_2$ then
\[
x_{r+2 \to r+1} x_{r+2 \leftarrow r+1} (e^2_{r+1}) = x_{r+1 \leftarrow
  r} x_{r+1 \to r} (e^2_{r+1}) = x_{r+1 \leftarrow r} (e^2_r) \ne 0.
\]
Therefore, $x_{r+2 \leftarrow r+1} (e^2_{r+1}) \ne 0$.  But $x_{r+2
\leftarrow r+1} (e^2_{r+1})$ must be a multiple of $e^1_{r+2}$ as
above.  Thus we must have $r+2 \le k_1$ and $x_{r+2 \leftarrow r+1}
(e^2_{r+1}) \ne 0$.  Continuing in
this manner we see that $k_2 < k_1$.
Refer to Figure~\ref{fig:youngcomm} for illustration.
\begin{figure}
\begin{center}
\epsfig{file=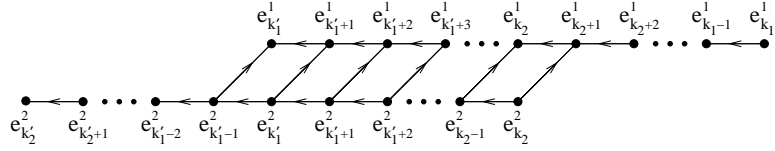,width=4in}
\caption{If $x_{r+1 \leftarrow r} (e^2_r) \ne 0$ for some $r$, the
  commutativity of
  this diagram forces $k_2' < k_1'$ and 
  $k_2 < k_1$. Vertices represent the spans of the indicated vectors.
  Those aligned vertically lie in the same
  $\mathbf{V}_i$. The arrows indicate the action of the obvious
  component of $x$.  \label{fig:youngcomm}}
\end{center}
\end{figure}

Now, let $x$ lie in the conormal bundle to the point
\begin{equation}
\bigoplus_{i=1}^s x(k_i',
  k_i'+l_i-1) \in \mathbf{E}_{\oplus_{i=1}^s
  \mathbf{V}_\infty(k_i', k_i'+l_i-1), \Omega}.
\end{equation}
We can assume (by reordering the indices if necessary) that $k_1' \ge
k_2' \ge \dots \ge k_s'$.  Now, by the above arguments, $x_{r+1
  \leftarrow r} (e^i_r)$ must be a linear combination of
$\{e^j_{r+1}\}_{j < i}$.  Thus
\[
e^1_{k_1'} \in \ker x_{k_1' \to k_1'-1} \cap \ker x_{k_1' + 1
  \leftarrow k_1'}. 
\]
By the stability condition, we must then have $k_1'=0$ and there can be
no other $e^i_r$ in $\ker x_{r \to r-1} \cap \ker x_{r+1 \leftarrow
r}$ for any $r$.  Now, by the above considerations, $e^2_{k_2'}$ is in
$\ker x_{k_2' \to k_2'-1} \cap \ker x_{k_2'+1 \leftarrow k_2'}$ unless
$k_2'+1=k_1'$ and $x_{k_1' \leftarrow k_2'} (e^2_{k_2'})$ is a
non-zero multiple of $e^1_{k_1'}$.  Continuing in this manner, we see
that we must have $k_{i+1}'+1 = k_i'$ and $x_{k_i' \leftarrow k_{i+1}'}
(e^{i+1}_{k_{i+1}'}) = c_i e^i_{k_i'} \ne 0$ for $1 \le i \le s-1$.
Then by the above we must have $k_{i+1} < k_i$ for $1 \le i \le s-1$.
Setting $l_i = k_i - k_i' +1$ the theorem follows.
\end{proof}

\begin{figure}
\begin{center}
\epsfig{file=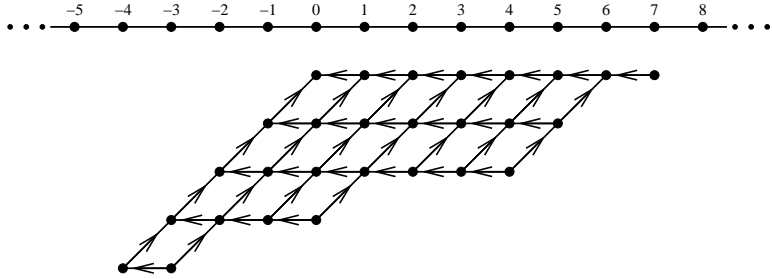,width=4in}
\caption{The irreducible components of $\mathcal{L}(\mathbf{v}, \mathbf{w}^0)$ are enumerated by Young diagrams.  The top line is the Dynkin graph of type
$A_\infty$.  The other horizontal lines represent $x_\infty(k',k)$ where $k'$
and $k$ are the positions of the leftmost and rightmost vertices.
\label{fig:youngdiag}}
\end{center}
\end{figure}

The Young diagrams enumerating the irreducible components of
$\mathcal{L}(\mathbf{v}, \mathbf{w}^0)$ can be visualized as in
Figure~\ref{fig:youngdiag}.  Note that the vertices in our diagram
correspond to the boxes in the classical Young diagram, and our arrows
intersect the classical diagram edges.

For the level one $A_\infty$ case, it is relatively easy to compute
the geometric action of the generators $E_k$ and $F_k$ of
$\mathfrak{g}$.  We first note that for every $\mathbf{v}$,
$\mathcal{L}(\mathbf{v},\mathbf{w}^0)$ is either empty or is a point.
It follows that each $X_Y$ is equal to
$\mathcal{L}(\mathbf{v},\mathbf{w}^0)$ for some unique $\mathbf{v}$
which we will denote $\mathbf{v}_Y$.

\begin{lem}
The function $g_{X_Y}$ corresponding to the irreducible component
$X_Y$ where $Y \in \mathcal{Y}$ is simply $\id_{X_Y}$, the function on
$X_Y$ with constant value one.
\end{lem}

\begin{proof}
This is obvious since $X_Y$ is a point.
\end{proof}

\begin{prop}
One has $F_k \id_{X_Y} = \id_{X_{Y'}}$ where $\mathbf{v}_{Y'} = \mathbf{v}_Y +
\mathbf{e}^k$ if such a $Y'$ exists and $F_k \id_{X_Y} = 0$
otherwise.  Also, $E_k \id_{X_Y} = \id_{X_{Y''}}$ where $\mathbf{v}_{Y''} =
\mathbf{v}_Y - \mathbf{e}^k$ if such a $Y''$ exists and $E_k \id_{X_Y}
= 0$ otherwise.
\end{prop}

\begin{proof}
It is clear from the definitions that $F_k \id_{X_Y} = c_1
\id_{X_{Y'}}$ and $E_k \id_{X_Y} = c_2 \id_{X_{Y''}}$ for some
constants $c_1$ and $c_2$ if $Y'$ and $Y''$ exist as described above
and that these actions are zero otherwise.  We simply have to compute
the constants $c_1$ and $c_2$. Now,
\begin{align*}
F_k \id_{X_Y} (x) &= (\pi_2)_! \pi_1^* \id_{X_Y}(x) \\
&= \chi (\{S\, |\, S \text{ is $x$-stable},\, x|_S \in X_Y\}) \\
&= \chi (\text{pt}) \\
&= 1
\end{align*}
if $x \in X_{Y'}$ where $\mathbf{v}_{Y'} = \mathbf{v}_Y +
\mathbf{e}^k$ and zero otherwise.  The fact that the above set is
simply a point follows from the fact that $S_k$ must be the sum of the
images of
$x_h$ such that $\inc(h)=k$.  Thus $c_1=1$ as desired.

Note that if there exists a $Y'$ such that $\mathbf{v}_{Y'} = \mathbf{v}_Y +
\mathbf{e}^k$ then there cannot exist a $Y''$ such that $\mathbf{v}_{Y''} =
\mathbf{v}_Y - \mathbf{e}^k$ and vice versa.  Therefore if such a
$Y''$ exists, $F_k \id_{X_Y} =0$ and so
\[
H_k \id_{X_Y} = [E_k,F_k] \id_{X_Y} = -F_k E_k \id_{X_Y}.
\]
One can easily check that $H_k \id_{X_Y} = -\id_{X_Y}$ if a $Y''$
exists as described above and thus $F_k E_k \id_{X_Y} = \id_{X_Y}$.
It then follows from the above that we must have $c_2=1$.
\end{proof}

The above action of the type $A_\infty$ Lie algebra in the space
spanned by a basis indexed by Young diagrams is well known in a purely
algebraic context (see e.g. \cite{JM}).

\begin{rem}{\upshape
All the results of this section can be repeated with minor
modifications for fundamental representations of finite-dimensional
Lie algebras of type $A_n$.  In this case, the bases of fundamental
representations will be enumerated by Young diagrams of size bounded
by an $m \times (n+1-m)$ rectangle, where $m=1,2,\dots,n$ is the
index of the fundamental representation.  Note that the same Young diagrams
also naturally enumerate the Schubert cells of the Grassmannians
$Gr(m, n+1)$ for type $A_n$ or the semi-infinite Grassmannian for type
$A_\infty$.}
\end{rem}


\subsection{Type $A_n^{(1)}$}
\label{sec:Anlevel1}
Let $\mathcal{Y}_n$ be the set of all Young diagrams $[l_1, \dots,
l_s]$ satisfying $l_i > l_{i+n}$ for all $i=1, \dots, s$ ($l_j=0$ for
$j>s$).  For $Y=[l_1, \dots, l_s] \in \mathcal{Y}_n$, let $A_Y$ be the
set $\{(1-i,l_i-i)\ |\ 1 \le i \le s\}$.

\begin{theo}
The irreducible components of $\mathcal{L}(\mathbf{v},\mathbf{w}^0)$
are precisely those $X_f$ where $f \in \tilde Z_\mathbf{V}$ such that
\[
\{(k',k)\ |\ f(k',k)=1\} = A_Y
\]
for some $Y \in \mathcal{Y}_n$ and $f(k',k)=0$ for $(k',k) \not \in
A_Y$ (up to simultaneous translation of $k'$ and $k$ by $n+1$).
Denote the component corresponding to such an $f$ by $X_Y$.  Thus,
$Y \leftrightarrow X_Y$ is a natural 1-1 correspondence between the set
$\mathcal{Y}_n$ and the
irreducible components of $\cup_\mathbf{v} \mathcal{L}(\mathbf{v},
\mathbf{w}^0)$.
\end{theo}

\begin{proof}
The argument is exactly analogous to that used in the proof of
Theorem~\ref{thm:irrcomp_lev1_infty}.  We need only note that a
point in the conormal bundle to the orbit through the point
\begin{equation}
\label{eq:sum_reps}
\sum_{i=1}^s x(k_i', k_i'+l_i-1) \in \mathbf{E}_{\oplus_{i=1}^s
  \mathbf{V}(k_i', k_i'+l_i-1), \Omega}
\end{equation}
lies in
$\Lambda_{\mathbf{V}}(\mathbf{v}, \mathbf{w}^0)$ if
and only if $l_i > l_{i+n}$ for all $i=1,\dots,s$ ($l_i =0$ for $i>s$)
by the aperiodicity condition.
\end{proof}

Note that Nakajima's construction yields an
action of the Lie algebra on the basis $\{g_{X_Y}\}_{Y \in
  \mathcal{Y}_n}$ of the basic representation.  However, this action
is not as straightforward to compute as in the $A_\infty$ case and will
be considered in a future work.


\subsection{$\widehat{\mathfrak{gl}}_{n+1}$ Case}
\label{sec:rem_gl_lev1}

We define $A_Y$ for $Y \in \mathcal{Y}$ as in Section~\ref{sec:Ainflevel1}.

\begin{theo}
The irreducible components of $\tilde{\mathcal{L}}(\mathbf{v},\mathbf{w}^0)$
are precisely those $X_f$ where $f \in \tilde Z_\mathbf{V}$ such that
\[
\{(k',k)\ |\ f(k',k)=1\} = A_Y
\]
for some $Y \in \mathcal{Y}$ and $f(k',k)=0$ for $(k',k) \not \in
A_Y$ (up to simultaneous translation of $k'$ and $k$ by $n+1$).
Denote the component corresponding to such an $f$ by $X_Y$.  Thus,
$Y \leftrightarrow X_Y$ is a natural 1-1 correspondence between the set
$\mathcal{Y}$ and the
irreducible components of $\cup_\mathbf{v} \tilde{\mathcal{L}}(\mathbf{v},
\mathbf{w}^0)$.
\end{theo}

\begin{proof}
The argument is exactly analogous to that used in the proof of
Theorem~\ref{thm:irrcomp_lev1_infty}.
\end{proof}

As noted in Section~\ref{sec:rem_gl_lus}, since
$\widehat{\mathfrak{gl}}_{n+1}$ is not a Kac-Moody algebra we need to
modify Nakajima's construction of highest weight representations.
Note that for any $n$, the difference between
$\widehat{\mathfrak{gl}}_{n+1}$ and $\widehat{\mathfrak{sl}}_{n+1}$ is
the same Heisenberg algebra $\widehat{\mathfrak{gl}}_1$.  The
representations of Heisenberg algebras in the context of geometric
representation theory were first constructed by Grojnowski \cite{G96}
and Nakajima \cite{N97} (see \cite{N99} for a review).  However, it is
not obvious how to adapt this representation theory to the new quiver
varieties $\tilde{\mathcal{L}}(\mathbf{v}, \mathbf{w}^0)$, obtaining
the desired commutation relations with the
generators of $\widehat{\mathfrak{sl}}_{n+1}$.  This problem will be
considered in a future work.


\section{Arbitrary Level Representations}
\label{sec:arblevel}
\setcounter{theo}{0}

\subsection{Type $A_\infty$}

We first recall some definitions from \cite{D892}.  A \emph{Maya
diagram} is a bijection $m: \Z \to \Z$ such that $(m(j))_{j < 0}$
and $(m(j))_{j \ge 0}$ are both increasing.  For each Maya diagram
there exists a unique $\gamma \in \Z$ such that $m(j) - j = \gamma$
for $|j| \gg 0$.  This $\gamma$ is called the \emph{charge} of $m$.
We denote the set of Maya diagrams of charge $\gamma$ by $\mathcal{M}[\gamma]$.
For $m \in \mathcal{M}[\gamma]$ we let
\[
m[r] = (m(j)+r)_{j \in \Z} \in \mathcal{M}[\gamma + r].
\]

\begin{figure}
\begin{center}
\epsfig{file=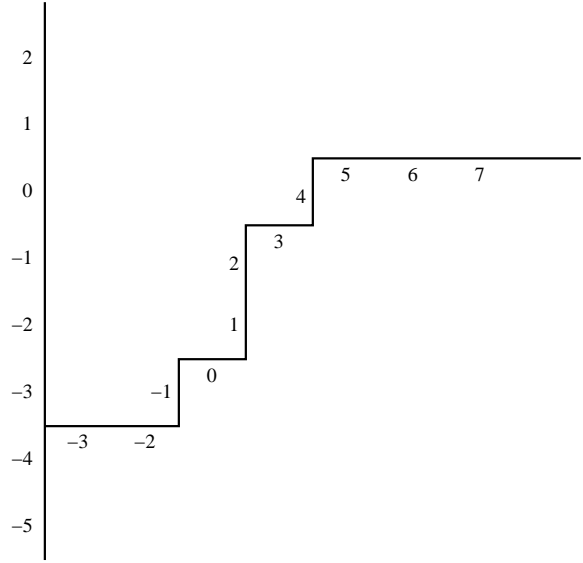,width=3in}
\caption{The Maya diagram corresponding to
$(m(j))_{j\ge0} = (-3,-2,0,3,5,6,7,\dots)$,
$(m(j))_{j<0} = (\dots,-6,-5,-4,,-1,1,2,4)$.
\label{fig:mayadiag_eg}}
\end{center}
\end{figure}

We can visualize a Maya diagram by a Young diagram.  Consider a
lattice on the right half plane with lattice points $\{(i,j) \in
\Z^2\, |\, i \ge 0\}$.  Each edge on the lattice is oriented, starting
at $(i,j)$ and ending at $(i+1,j)$ or $(i,j+1)$ and is numbered by the
integer $i+j$.  A \emph{path} on the lattice is a map $e$ from $\Z$ to
the set of edges on the lattice such that $e(j)$ has number $j$ and
the ending site of $e(j)$ is the starting site of $e(j+1)$.  To each
Maya diagram of charge $\gamma$, we associate the unique path
satisfying the following conditions.
\begin{enumerate}
\item For $j \ll 0$, $e(j)$ is the edge from $(0,j)$ to $(0,j+1)$,

\item The edge $e(m(j))$ is vertical (resp. horizontal) if $j<0$
  (resp. $j \ge 0$).
\end{enumerate}
Note that these conditions imply that for $j \gg 0$, $e(j)$ is the
edge from $(j- \gamma,\gamma)$ to $(j-\gamma+1,\gamma)$.
See Figure~\ref{fig:mayadiag_eg}.

Such a path divides the right half plane into two components.  The
upper half is an \emph{infinite Young diagram} $\mathfrak{Y}$ which
consists of a quadrant and a (finite) Young diagram $Y$ attached along
a horizontal line at height $\gamma$.  Thus the set of Maya diagrams
are in one-to-one correspondence with the set of pairs $(Y,\gamma)$
where $Y \in \mathcal{Y}$ and $\gamma \in \Z$.

\begin{lem}[\cite{D892}]
\label{lem:maya_order}
Let $m \in \mathcal{M}[\gamma]$, $m' \in \mathcal{M}[\gamma']$, and
let $\mathfrak{Y}$, $\mathfrak{Y}'$ be the corresponding infinite
Young diagrams.  Then the following conditions are equivalent.
\begin{enumerate}
\item $m(j) \le m'(j)$ for $j \ge 0$,
\item $\gamma \le \gamma'$ and $m(j - \gamma) \ge m'(j- \gamma')$ for
  $j < \gamma$,
\item $\mathfrak{Y} \supset \mathfrak{Y}'$.
\end{enumerate}
\end{lem}

We put a partial ordering on the set of Maya diagrams by letting $m
\le m'$ if the conditions in Lemma~\ref{lem:maya_order} hold.

Let $\Lambda = \Lambda_{\gamma_1} + \dots
\Lambda_{\gamma_l}$ where
$\gamma_1 \le \dots \le \gamma_l$ and the $\Lambda_i$ are
fundamental weights of $\mathfrak{g}$.  Let $\mathbf{w} \in (\Z_{\ge
  0})^\Z$ (that is, $\mathbf{w}$ is function from $\Z$ to $\Z_{\ge 0}$) be
the vector with $i^{th}$ component equal to the number of $\gamma_j$
equal to $i$.  Let
\[
\mathcal{M}[\Lambda] = \mathcal{M}[\gamma_1] \times \dots \times
\mathcal{M}[\gamma_l].
\]
For $Y=[l_1,\dots,l_s] \in \mathcal{Y}$, let $A^\gamma_Y$ be the set
$\{(\gamma+1-i, \gamma + l_i-i)\, |\, 1 \le i \le s\}$.
For $M = ((Y_1,\gamma_1), \dots, (Y_l,\gamma_l)) \in
\mathcal{M}[\Lambda]$, let $A_M = \cup_{i=1}^l A^{\gamma_i}_{Y_i}$ and let
$f_M \in \tilde{Z}^\infty$ be the function such that $f(k',k)$ is
equal to the number of
times $(k',k)$ appears in the set $A_M$.

\begin{theo}
\label{thm:irrcomp_arb_infty}
The irreducible components of $\mathcal{L}(\mathbf{v},\mathbf{w})$
are precisely those $X_f$ where $f = f_M$
for some $M \in \mathcal{M}[\Lambda]$.
Denote the component $X_{f_M}$ by $X_M$.
Then $M \leftrightarrow X_M$ is a natural 1-1 correspondence between the set
\[
\{(m_1, \dots, m_l) \in \mathcal{M}[\Lambda]\, |\, m_1 \le \dots \le
m_l \}
\]
and the
irreducible components of $\cup_\mathbf{v} \mathcal{L}(\mathbf{v},
\mathbf{w})$.
\end{theo}

\begin{proof}
Recall that irreducible components of
$\mathcal{L}(\mathbf{v},\mathbf{w})$ are the closures of the
$G_\mathbf{V}$-orbits (or isomorphism classes) in
$\mathbf{E}_{\mathbf{V}, \Omega}$ and that there is a representative
of each orbit of the form
\begin{equation}
\label{eq:xdecomp}
x = \bigoplus_{(k' \le k) \in K} x_\infty(k',k)
\end{equation}
for some finite set of pairs $K$.
By picturing $x_\infty(k',k)$ as the string of vertices $k', k'+1,
\dots, k$, we can represent such an $x$ by a set of finite strings of
vertices corresponding to the various $x_\infty(k',k)$ appearing in
\eqref{eq:xdecomp}.  We call the number of
vertices in a string its \emph{length}.  Each vertex of a
string represents a basis vector of $\mathbf{V}$ with degree given by
the location of the vertex.  The action of $x$ maps each of these
basis vectors to the basis vector corresponding to the next (one
lower) vertex in
the string (or to zero if no such vertex exists).
See Figure~\ref{fig:strings_eg}.

\begin{figure}
\begin{center}
\epsfig{file=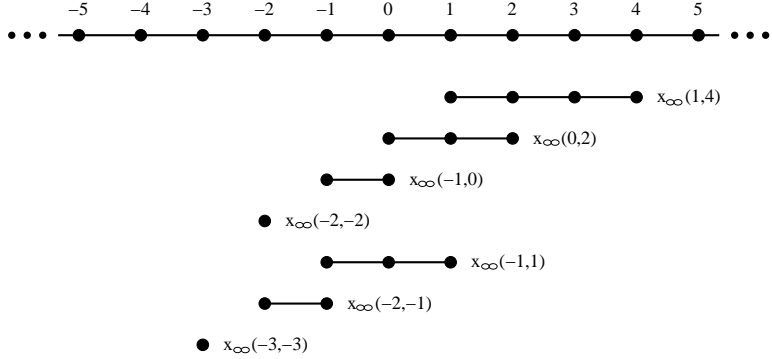,width=4in}
\caption{The strings associated to some $x \in \mathbf{E}_{\mathbf{V},
\Omega}$. The top line is the Dynkin diagram of type $A_\infty$.
\label{fig:strings_eg}}
\end{center}
\end{figure}

It is then a straightforward extension of the proof of
Theorem~\ref{thm:irrcomp_lev1_infty} that the allowable sets of
strings are precisely those that can be grouped into subsets, one for
each $\gamma_i$, such that the subset corresponding to $\gamma_i$,
when ordered by decreasing leftmost vertex, has weakly decreasing
lengths, the first leftmost vertex is $\gamma_i$ and the leftmost
vertices decrease by one as we move through the subset in order
(by leftmost, we mean the vertex with the smallest index).  This
is precisely the first claim of the Theorem.

It is easy to see that  many different $M \in
\mathcal{M}[\Lambda]$ may correspond to the same irreducible
component.  For example, for $\Lambda = \Lambda_{-1} + \Lambda_1$, both
\begin{gather*}
M = (([3,2,1],-1),([4,3,2,1],1)),\ \text{and} \\
M' = (([2,1],-1),([4,3,3,2,1],1))
\end{gather*}
belong to $\mathcal{M}[\Lambda]$ and correspond to the set of strings
shown in Figure~\ref{fig:strings_eg} (and hence to the same
irreducible component).
However, we can associate a unique $M \in \mathcal{M}[\Lambda]$ to
each set of strings described above as
follows.  Associate to $\gamma_1$ the longest string with leftmost
vertex $\gamma_1$ and remove this string from the set.  Now do the
same for $\gamma_2$, etc.  After we have associated a string to
$\gamma_l$, we start again at $\gamma_1$ but this time we select the
longest string with leftmost vertex $\gamma_1-1$ and so forth.  If at
any point, there is no string to associate with a given $\gamma_i$, we
remove this $\gamma_i$ from further steps.  In this way we associate
to each $\gamma_i$ a sequence of strings of weakly decreasing length
(by our condition on the possible sets of strings) with leftmost
vertices decreasing by one.  The lengths of the strings associated to
$\gamma_i$ give a Young diagram $Y_i$ and we set $m_i=(Y_i,\gamma_i)$.
By construction, the length of any string
associated to $\gamma_i$ is greater than the length of a string with
the same left end point associated to $\gamma_j$ for $j>i$.  This
immediately yields the condition $m_1 \le \dots \le m_l$.  Our
construction thus gives us the one-to-one correspondence asserted in
the Theorem.
\end{proof}

Note that the enumeration of the irreducible components given in
Theorem~\ref{thm:irrcomp_arb_infty} matches that of Proposition~4.6 of
\cite{D892}.


\subsection{Type $A_n^{(1)}$}
\label{sec:arblevel_an1}
We now consider the case where $\mathfrak{g}$ is of type $A_n^{(1)}$.
For an element $M = (m_1, \dots, m_l) \in
\mathcal{M}[\Lambda]$, let $R_M$ be the set (with multiplicity) of
pairs $(i,l_i)$ where
$l_i$ is the length of a row with top edge having $y$-coordinate $i$
belonging to one of the $m_j$.  We say that $M$ is \emph{$n$-reduced}
if
\[
\{(k+i,l)\, |\, 0 \le i \le n\} \not \subset R_M
\]
for all $k$ and $l$.

Define $f_M$ for $M \in \mathcal{M}[\Lambda]$ as in the previous subsection
(except that now our pairs are defined only up to simultaneous translation
by $n+1$).

\begin{theo}
\label{thm:irrcomp_an}
The irreducible components of $\mathcal{L}(\mathbf{v},\mathbf{w})$
are precisely those $X_f$ where $f = f_M$
for some $n$-reduced $M \in \mathcal{M}[\Lambda]$.
Denote the component $X_{f_M}$ by $X_M$.  Then
$M \leftrightarrow X_M$ is a natural one-to-one correspondence between the set
\[
\{(m_1, \dots, m_l) \in \mathcal{M}[\Lambda]\, |\, m_1 \le \dots \le
m_l \le m_1[n+1],\, \text{$M$ is $n$-reduced} \}
\]
and the
irreducible components of $\cup_\mathbf{v} \mathcal{L}(\mathbf{v},
\mathbf{w})$.
\end{theo}
\begin{proof}
Each irreducible component corresponds to a set of strings as in the
proof of Theorem~\ref{thm:irrcomp_arb_infty} with the added condition
that we cannot have $n+1$ strings, each of the same length, whose left
end points are the $n+1$ vertices of our quiver.  That is, we must have
that $M$ is $n$-reduced.  Note that the process
described in the proof of Theorem~\ref{thm:irrcomp_arb_infty} yields
$m_i=(Y_i,\gamma_i)$
satisfying $m_1 \le \dots \le m_l \le m_1[n+1]$ as desired.  The Theorem
follows.
\end{proof}

Again, as noted in Section~\ref{sec:Anlevel1}, Nakajima's construction
yields an
action of the Lie algebra on the bases $\{g_{X_M}\}$ of the
irreducible representations in both the $A_\infty$ and $A_n^{(1)}$
cases which is more difficult to directly compute than in the level 1
$A_\infty$ case.


\subsection{$\widehat{\mathfrak{gl}}_{n+1}$ Case}

\begin{theo}
\label{thm:irrcomp_gn}
The irreducible components of $\tilde{\mathcal{L}}(\mathbf{v},\mathbf{w})$
are precisely those $X_f$ where $f = f_M$
for some $M \in \mathcal{M}[\Lambda]$.
Denote the component $X_{f_M}$ by $X_M$.  Then
$M \leftrightarrow X_M$ is a natural one-to-one correspondence between the set
\[
\{(m_1, \dots, m_l) \in \mathcal{M}[\Lambda]\, |\, m_1 \le \dots \le
m_l \le m_1[n+1] \}
\]
and the
irreducible components of $\cup_\mathbf{v} \tilde{\mathcal{L}}(\mathbf{v},
\mathbf{w})$.
\end{theo}
\begin{proof}
The argument is the same as the proof of Theorem~\ref{thm:irrcomp_an}
except that we do not have the aperiodicity condition and thus do not
require that $M$ is $n$-reduced.
\end{proof}

Note that the enumeration of the irreducible components of
$\cup_\mathbf{v} \tilde{\mathcal{L}}(\mathbf{v},\mathbf{w})$ given by
Theorem~\ref{thm:irrcomp_gn} is the same as that given
by Proposition~4.7 of \cite{D892} for a
spanning set of the dual to the irreducible highest weight
representation of $\widehat{\mathfrak{gl}}_{n+1}$.
In order to extend the geometric construction of highest weight
representations of $\widehat{\mathfrak{sl}}_{n+1}$ to
$\widehat{\mathfrak{gl}}_{n+1}$ for an arbitrary level, one would need a
representation of the Heisenberg algebra as discussed in
Section~\ref{sec:rem_gl_lev1}.  Here one might use the construction of
the Heisenberg algebra by Baranovsky \cite{B} that generalizes the
Grojnowski/Nakajima construction to higher levels.

\begin{rem}{\upshape
One can also give a geometric interpretation of the full Fock space of
\cite{D892} with basis indexed by $\mathcal{M}[\Lambda]$ via the
``smooth'' $U_l$-instanton moduli space $\sqcup_r \mathcal{M}(r,l)$
which has the same generating function for cohomology (see
e.g. \cite{N99}, Chapter 5) as the full Fock space with the basis
$\mathcal{M}[\Lambda]$.  The types $A_n^{(1)}$ or $A_\infty$ are
reflected in the respective action of the groups $\Z/(n+1)\Z$ or $\C^*$
on the moduli space, and $\gamma_1, \dots, \gamma_l$ is the set of
one-dimensional representations of these groups that determine this
action.}
\end{rem}


\section{A Comparison With The Path Space Representation}
\label{sec:comparison}
\setcounter{theo}{0}

The authors of \cite{D89} constructed the basic representation of
$A_n^{(1)}$ on the space of paths.  In \cite{D892}, this path
realization is generalized to arbitrary level.
We now compare the geometric presentation $L(\mathbf{v},\mathbf{w})$
with theirs.  We will slightly modify the definitions of \cite{D89} to
agree with the more general definitions of \cite{D892}.


\subsection{The Level One Case}
\label{sec:comp_level1}

A \emph{basic path} is a sequence $p=(\lambda_0, \lambda_1, \dots)$ of
integers $\lambda_i \in \{0,1,\dots,n\}$.  The basic path
\[
(\overline{j})_{j \ge 0} = (0,1,\dots,n,0,1,\dots,n,\dots)
\]
is called the \emph{ground state}.  Here $\overline{k}$
for $k \in \Z$ signifies the unique integer such that $0 \le
\overline{k} \le n$ and $\overline k = k$ mod $n+1$.
Let
\[
\mathcal{P}_b = \{p=(\lambda_0, \lambda_1,\dots)\ |\ \lambda_j =
\overline{j} \text{ for all but a finite number of $j$}\}.
\]

For a basic path $p = (\lambda_0, \lambda_1, \dots) \in
\mathcal{P}_b$, let
\[
\omega(p) = \sum_{i=1}^\infty i(H(\lambda_i, \lambda_{i+1}) -
H(\overline{i}, \overline{i+1})),
\]
where
\[
H(\lambda, \mu) =
  \begin{cases}
    0 & \text{if $\lambda < \mu$} \\
    1 & \text{if $\lambda \ge \mu$}
  \end{cases}.
\]

Basic paths in $\mathcal{P}_b$ can be labeled by Young diagrams as we now
describe.  The set $\mathcal{M}[0]$ is in one-to-one correspondence
with the set of strictly increasing sequences 
of integers $m=(m(0), m(1), \dots)$ such that $m(j) = j$ for $j$ large
and
\[
\#\{j\ |\ m(j) <0\} = \#(\{0,1,2,\dots\} - \{m(j)\ |\ m(j) \ge 0\}).
\]
Such a sequence represents the Young diagram of signature $[\dots
3^{r_3} 2^{r_2} 1^{r_1}]$ where $r_j = m(j) - m(j-1) -1$ and
vice versa.
To a Maya diagram $m=(m(0), m(1), \dots) \in \mathcal{Y}_n$
we associate the basic path $p=(\overline{m(0)}, \overline{m(1)},
\dots) \in \mathcal{P}_b$.  Then the ground state corresponds to the
empty Young diagram $\phi$.  In the sequel, we will identify
$n$-reduced Young
diagrams (that is, elements of $\mathcal{Y}_n$) and basic paths via
the above correspondence.

For $Y = [l_1, \dots, l_s] \in \mathcal{Y}_n$, let
\[
\Delta_k(Y) = \delta(k,-s) + \sum_{i=1}^s \left( \delta(k,l_i -i+1) -
  \delta(k, l_i-i) \right).
\]

\begin{prop}
One has $H_k g_{X_Y} = \Delta_k(Y) g_{X_Y}$.
\end{prop}
\begin{proof}
Let $Y=[l_1, \dots, l_s]$.  Then $g_{X_Y} \in
L(\mathbf{v},\mathbf{w}^0)$ where
\begin{align*}
\mathbf{v} &= \dim \bigoplus_{i=1}^s \mathbf{V}(1-i, l_i-i) \\
&= \sum_{i=1}^s \sum_{l=1-i}^{l_i-i} \alpha_{\bar l}.
\end{align*}
Recall that the weight of the space $L(\mathbf{v},\mathbf{w}^0)$ is
\[
(u_0, \dots, u_n) = \mathbf{w}^0 - C\mathbf{v}
\]
and thus $H_k g_{X_Y} = u_k g_{X_Y}$ with
\begin{align*}
u_k &= \Lambda_0(\alpha_k) - \sum_{i=1}^s \sum_{l=1-i}^{l_i-i} \left<
  \alpha_k, \alpha_{\bar l} \right> \\
&= \delta(k,0) - \sum_{i=1}^s \sum_{l=1-i}^{l_i-i} \left( 2\delta(k,l)
  - \delta(k,l-1) - \delta(k,l+1) \right) \\
&= \delta(k,0) - \sum_{i=1}^s \left( \delta(k,1-i) - \delta(k,-i) +
  \delta(k,l_i-i) - \delta(k,l_i-i+1) \right) \\
&= \delta(k,-s) + \sum_{i=1}^s \left( \delta(k,l_i-i+1) -
  \delta(k,l_i-i) \right) \\
&= \Delta_k(Y).
\end{align*}
\end{proof}

\begin{prop}
One has $d(g_{X_Y}) = -\omega(Y) g_{X_Y}$.
\end{prop}

\begin{proof}
We first compute the left hand side.  It is obvious that
\[
d(g_{X_Y}) = -\mathbf{v}_0 g_{X_Y}
\]
where $X_Y \subset \mathcal{L}(\mathbf{v}, \mathbf{w}^0)$.
Consider the representation $(\mathbf{V}(k',k'+l-1),x(k',k'+l-1))$ where
$l=(n+1)a +b$ with $0 \le b \le n$.  Then
\begin{align*}
\mathbf{v}_0 &= \dim \mathbf{V}(k',k'+l-1)_0 \\
&= a + \begin{cases}
  1 & \text{ if } \overline{k'-1}+b>n \\
  0 & \text{ if } \overline{k'-1}+b \le n
\end{cases}.
\end{align*}

Thus, for $Y=[l_1,\dots,l_s] \in \mathcal{Y}_n$ where $l_i = (n+1)a_i
+ b_i$ with $0 \le b_i \le n$,
\[
\mathbf{v}_0 = \sum_{i=1}^s \left( a_i +
\begin{cases}
  1 & \text{if } \overline{-i} + b_i > n \\
  0 & \text{otherwise}
\end{cases}
\right).
\]

We now show that this is equal to $\omega(Y)$.  Let $Y_i=[l_1, \dots,
l_i]$ for $0 \le i \le s$ where $Y_0 = \phi$ and let $(\lambda^i_0,
\lambda^i_1, \dots)$ be the corresponding basic path.  Then the first $l_i$
positions of the basic path corresponding
to $Y_{i-1}$ are 
\[
(\overline{1-i}, \overline{1-i}+1, \dots,
n,0,1,\dots,n,0,1,\dots,n,\dots,0,1,\dots,n,
0,1,\dots,\overline{b_i-i}).
\]
Here there are $a_i$ repetitions of $0,1,\dots,n$ if $\overline{i-1} <
b_i$ and $a_i-1$ repetitions if $\overline{i-1} \ge b_i$.

The first $l_i$ positions of the basic path corresponding to $Y_i$ are
simply obtained from the above by lowering all the entries by 1
(interpreting -1 as $n$).  The entries of $Y_i$ and $Y_{i-1}$ numbered
$l_i+1$ and above are
equal.  Then by considering the cases $\overline{i-1} <
b_i$ and $\overline{i-1} \ge b_i$, we see that
\[
\sum_{j=1}^\infty j(H(\lambda^i_j, \lambda^i_{j+1}) -
H(\lambda^{i-1}_j, \lambda^{i-1}_{j+1})) = a_i +
\begin{cases}
  1 & \text{if } \overline{-i} + b_i > n \\
  0 & \text{otherwise}
\end{cases}.
\]
and the result follows.
\end{proof}

\begin{theo}
The map $g_{X_Y} \mapsto Y$ is a weight-preserving vector space isomorphism
between the geometric presentation $L(\mathbf{w}^0)$ of
$L(\Lambda_0)$ and the basic path space representation given in \cite{D89}. 
\end{theo}

\begin{proof}
This follows directly from the previous two propositions and the
action of the $H_i$ and $d$ given in \cite{D89}.
\end{proof}


\subsection{Arbitrary Level}

We first recall some definitions of \cite{D892}.  
Let $\epsilon_\mu = (0,\dots,\overset{\mu-\text{th}}{1},\dots,0)$
for $0 \le \mu \le n$ denote the standard basis vectors of
$\Z^{n+1}$.  We extend the indices to $\Z$ by setting
$\epsilon_{\mu + n+1} = \epsilon_\mu$.  Fix a positive integer $l$
(the level of our representation).
A \emph{path} is a
sequence $\eta=(\eta(k))_{k\ge 0}$ consisting of elements $\eta(k) \in
\Z^{n+1}$ of the form $\epsilon_{\mu_1(k)} + \dots +
\epsilon_{\mu_l(k)}$ with $\mu_1(k), \dots, \mu_l(k) \in \Z$.
To a level $l$ dominant integral weight $\Lambda = \Lambda_{\gamma_1}
+ \dots \Lambda_{\gamma_l}$ is associated the path
\[
\eta_\Lambda = (\eta_\Lambda(k))_{k \ge 0},\ \eta_\Lambda(k) =
\epsilon_{\gamma_1+k} + \dots + \epsilon_{\gamma_l+k}.
\]
We call $\eta$ a $\Lambda$-path if $\eta(k) = \eta_\Lambda(k)$ for
$k \gg 0$.  The set of $\Lambda$-paths is denoted by
$\mathcal{P}(\Lambda)$.  Define the weight $\lambda_\eta$ of $\eta$ by
\[
\lambda_\eta = \Lambda - \sum_{k \ge 0} \pi (\eta(k) -
\eta_\Lambda(k)) - \omega(\eta)\delta
\]
where
\[
\omega(\eta) = \sum_{k \ge 1} k \left( H(\eta(k-1),\eta(k)) -
  H(\eta_\Lambda(k-1), \eta_\Lambda(k)) \right).
\]
Here $\delta$ is the null root and $\pi$ is the $\Z$-linear map from
$\Z^{n+1}$ to the weight lattice of the Lie algebra of type
$A_n^{(1)}$ such that
$\pi(\epsilon_\mu) = \Lambda_{\mu+1} - \Lambda_\mu$ (here
$\Lambda_{n+1} = \Lambda_0$).  The function $H$ is defined as follows:
if $\alpha = \epsilon_{\mu_1} + \dots + \epsilon_{\mu_l}$ and
$\beta = \epsilon_{\nu_1} + \dots + \epsilon_{\nu_l}$ ($0 \le
\mu_i, \nu_i \le n$), then
\[
H(\alpha,\beta) = \min_{\sigma \in S_l}\sum_{i=1}^l \theta(\mu_i - \nu_{\sigma(i)})
\]
where $S_l$ is the permutation group on $l$ letters, and
\[
\theta(\mu) = \begin{cases}
1 & \text{if $\mu \ge 0$} \\
0 & \text{otherwise}
\end{cases}.
\]

Note that we have redefined the notation $\omega$ and $H$ of
subsection~\ref{sec:comp_level1}.  However, our new definitions reduce
to the old ones in the case $\Lambda=\Lambda_0$ and so, to avoid a
proliferation
of notation, we denote the new functions by the same symbols.

Let $\eta$ be a $\Lambda$-path.  An element $M = (m_1,\dots,m_l) \in
\mathcal{M}[\Lambda]$ is called a \emph{lift} of $\eta$ if and only if
\begin{equation}
\label{eq:liftdef1}
m_1 \le \dots \le m_l \le m_1[r]
\end{equation}
and
\begin{equation}
\label{eq:liftdef2}
\eta(k) = \epsilon_{m_1(k)} + \dots + \epsilon_{m_l(k)}.
\end{equation}

If $M=(m_1,\dots,m_l)$ and $M'=(m_1',\dots,m_l')$ are lifts of a
$\Lambda$-path $\eta$ then we say $M \ge M'$ if and only if $m_j \ge
m_j'$ for $1 \le j \le l$.

Recall the definition of $R_M$ given in
Section~\ref{sec:arblevel_an1}.  For $M, M' \in \mathcal{M}[\Lambda]$,
we say that $M$ is an \emph{$n$-reduction} of $M'$ if $R_M$ is
obtained from $R_{M'}$ by the removal of sets of the form $\{(k+i,l)\,
|\, 0 \le i \le n\}$ for some $k$ and $l$.

\begin{prop}
Suppose $M=(m_1,\dots,m_l)$ is an $n$-reduction of $M'=(m_1',\dots,m_l')$
and $m_1 \le \dots \le m_l \le m_1[n+1]$, $m_1' \le \dots \le m_l' \le
m_1'[n+1]$.  Then $M$ and $M'$ are lifts of the
same path and $M \ge M'$.
\end{prop}
\begin{proof}
Recall the construction in the proof of
Theorem~\ref{thm:irrcomp_arb_infty}.  Note that choosing arbitrary
strings instead of the longest string at each step will not change the
values of the right hand side of \eqref{eq:liftdef2} (for any $k$).
Thus, let us form $M''=(m_1'',\dots,m_l'') \in \mathcal{M}[\Lambda]$
from the same strings comprising $M'$ but where one of the $m_i''$
contains the entire set of strings of the form $\{(k+i,l)\,
|\, 0 \le i \le n\}$ which is removed from $R_{M'}$ to obtain $R_M$.
Now, removing this set of strings from $M''$ simply amounts to
removing this set from $m_i''$.  But this just cuts an $(n+1) \times
l$ square out of the Maya diagram $m_i''$ and shifts the part of the
diagram below the cut up $n+1$ units.  See Figure~\ref{fig:nreduction}.

\begin{figure}
\begin{center}
\epsfig{file=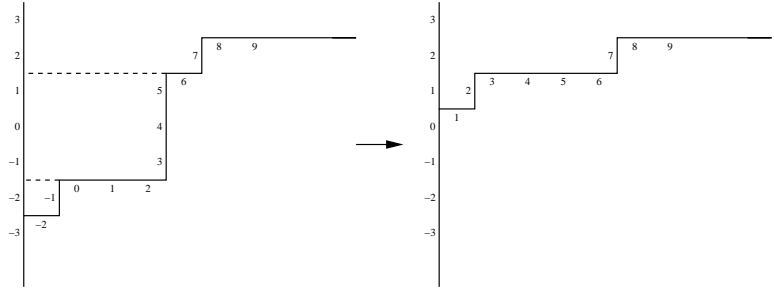,width=4in}
\caption{Removing an $(n+1) \times l$ square from a Maya diagram (here $n=2$ and $l=4$).  Notice that the enumeration of the horizontal edges does not
change mod $(n+1)$.
\label{fig:nreduction}}
\end{center}
\end{figure}

Since $\epsilon_{\mu+n+1} = \epsilon_\mu$, the values of the right
hand sides of
\eqref{eq:liftdef2} for $M'$ and $M''$ are the same.  However, $M$ is
simply obtained from $M''$ by applying the procedure of
Theorem~\ref{thm:irrcomp_arb_infty} to the strings of $M''$ and as
mentioned above, this does not change the value of the right hand sides of
\eqref{eq:liftdef2}.  Thus $M$ and $M'$ are lifts of the same path.

To show that $M \ge M'$, note that by the construction in the proof of
Theorem~\ref{thm:irrcomp_arb_infty}, $M$ is uniquely determined by
$R_M$.  Now, we obtain $R_M$ from $R_{M'}$ by
removing a set of the form $\{(k+i,l)\,
|\, 0 \le i \le n\}$ for some $k$ and $l$.  Thus, at each
stage in our construction of $M$, we chose a string of length less than
or equal to the string chosen in the construction of $M'$.  Thus we
have that $M \ge M'$.
\end{proof}

\begin{prop}[\cite{D892}]
\label{prop:highpathlifts}
For each $\Lambda$-path $\eta$ there exists a unique highest lift $M$
of $\eta$ such that $M \ge M'$ for any lift $M'$ of $\eta$.
\end{prop}

\begin{cor}
\label{cor:highpathlifts}
The set
\[
\{M=(m_1,\dots,m_l) \in \mathcal{M}[\Lambda]\, |\, M \text{ is
  $n$-reduced},\, m_1 \le \dots \le m_1 \le m_1[n+1]\} 
\]
is precisely the
set of highest lifts of paths in $\mathcal{P}(\Lambda)$.
\end{cor}

Let $M_\eta$ be the $n$-reduced element of $\mathcal{M}[\Lambda]$
corresponding to $\eta \in \mathcal{P}(\Lambda)$ and let
$\mathfrak{g}$ be the affine Lie algebra of type $A_n^{(1)}$.
Define
\[
\mathcal{P}(\Lambda)_\mu = \{\eta \in \mathcal{M}[\Lambda]\, |\,
\lambda_\eta = \mu\}.
\]
In \cite{D892}, the authors introduced a basis $\{\xi_\eta\, |\, \eta
\in \mathcal{P}(\Lambda)_\mu\}$ of the $\mu$ weight space of the
restricted dual of the
highest weight representation of $\mathfrak{g}$ of highest weight
$\Lambda$ \cite[Thm 5.4]{D892}.  The weight of $\xi_\eta$ is
$\lambda_\eta$ \cite[Thm 5.7]{D892}.

\begin{theo}
The map $g_{X_{M_\eta}} \mapsto \xi_\eta$ is a weight-preserving vector space
isomorphism between the geometric presentation
$L(\mathbf{w})$ of $L(\Lambda)$ and the path space representation of
\cite{D892}.
\end{theo}
\begin{proof}
The fact that we have a vector space isomorphism follows from
Proposition~\ref{prop:highpathlifts} and
Corollary~\ref{cor:highpathlifts}.  It remains to show that the map is
weight-preserving.
The definition of a path agrees with the definition of a basic path
when $\Lambda = \Lambda_0$ and the weights are the same in this case.
Thus we have the result for $\Lambda = \Lambda_0$ from the previous
subsection.  Then the result for arbitrary level one representations
follows easily.

Now, if
\[
M_\eta = ((Y_1,\gamma_1),\dots,(Y_l,\gamma_l))
\]
and $\mathbf{V}^i$ is the space corresponding to the strings (see the proof
of Theorem~\ref{thm:irrcomp_arb_infty}) of
$(Y_i,\gamma_i)$ (that is, its dimension in degree $j$ is equal to the
number of vertices of these strings that are numbered $j$) then the
weight of $g_{M_\eta}$ is
\[
\sum_{i=1}^l (\Lambda_{\gamma_i} - \dim \mathbf{V}^i)
\]
where $\dim \mathbf{V}^i$ is identified with an element of the root lattice as
in Section~\ref{sec:lus_def}.
But this is equal to $\sum_{i=1}^l
\lambda_{\eta_i}$ where $(Y_i,\gamma_i)$ is a lift of $\eta_i$ by the
level 1 result.  By
Proposition~5.6 of \cite{D892}, this is $\lambda_\eta$ as desired.
\end{proof}


\small{
}
\vspace{4mm}
\noindent
Igor B. Frenkel, \\
Department of Mathematics, Yale University, P.O. Box 208283, NEW HAVEN, CT, USA 06520-8283;\\
email:\;\texttt{frenkel-igor@yale.edu} \\
Alistair Savage,\\
Department of Mathematics, Yale University, P.O. Box 208283, NEW HAVEN, CT, USA 06520-8283;\\
email:\;\texttt{alistair.savage@aya.yale.edu}

\end{document}